\documentclass[reqno]{amsart}
\usepackage{amssymb,url}
\usepackage{hyperref}
\theoremstyle{plain}
\newtheorem{theorem}{Theorem}
\newtheorem{lemma}{Lemma}

\theoremstyle{remark}
\newtheorem{remark}{Remark}

 \allowdisplaybreaks[4]

\begin{document}

\title[A unified generalization of some quadrature rules]{A unified generalization of some quadrature rules and error bounds}

\author[W. J. Liu]{Wenjun Liu}
\address[W. J. Liu]{College of Mathematics and Physics\\
Nanjing University of Information Science and Technology \\
Nanjing 210044, China} \email{\href{mailto: W. J. Liu
<wjliu@nuist.edu.cn>}{wjliu@nuist.edu.cn}}

\author[J. Jiang]{Yong Jiang}
\address[J. Jiang]{College of Mathematics and Physics\\
Nanjing University of Information Science and Technology \\
Nanjing 210044, China}

\author[A. Tuna]{Adnan Tuna}
\address[A. Tuna]{Department of Mathematics\\
Faculty of Science and Arts\\
University of Ni\u{g}de\\
Merkez 51240, Ni\u{g}de, Turkey}
\email{atuna@nigde.edu.tr}

\begin{abstract}By introducing a parameter, we give a unified generalization of some quadrature rules, which not only unify the recent results about error bounds for generalized mid-point,
trapezoid and Simpson's rules, but also give some
new error bounds for other quadrature
rules as special cases. Especially, two sharp error inequalities are derived when $n$ is an odd and an even integer,
respectively.
\end{abstract}

\thanks{This paper was typeset using \AmS-\LaTeX}



\maketitle
\section{Introduction }

Error analysis for known and new quadrature rules has been
extensively studied in recent years. The approach from an
inequalities point of view to estimate the error terms has been used
in these studies (see \cite{bd}-\cite{u2006} and the references
therein).

In \cite{u2006}, by appropriately choosing the Peano kernel
\begin{equation}\label{1.1}
 T_n(x)=\left\{ \begin{array}{ll}\displaystyle
 \frac{(x-a)^{n-1}}{n!}\left[x+\frac{(n-2)a-nb}{2}\right],\ \ &\displaystyle x\in \left[a, \frac{a+b}{2}\right], \hfill
 \medskip\\
 \displaystyle \frac{(x-b)^{n-1}}{n!}\left[x+\frac{(n-2)b-na}{2}\right],\ \ &\displaystyle x\in \left(\frac{a+b}{2}, b\right], \hfill
\end{array} \right.
\end{equation}
error inequalities for a generalized trapezoid rule were given as follows.

\begin{theorem} \label{th1.1}
Let $f: [a, b] \rightarrow \mathbb{R}$ be a function such that $f^{(n-1)}$, $n > 1$, is absolutely continuous.

If
there exist real numbers $\gamma_n, \Gamma_n$ such that $\gamma_n\le f^{(n)}(x) \le \Gamma_n$, $x\in [a, b]$, then
\begin{align}\label{1.2}
 &\left|\int_a^b f(x)dx-\frac{f(a)+f(b)}{2}(b-a)+\sum_{i=1}^{ \lfloor \frac{n-1}{2} \rfloor}\frac{2i(b-a)^{2i+1}}{2^{2i}(2i+1)!}f^{(2i)}\left(\frac{a+b}{2}\right)\right|\nonumber\\
 \leq & \frac{\Gamma_n-\gamma_n}{(n+1)!}\frac{n}{2^{n+1}}(b-a)^{n+1},\quad
\text{if}\ n\ \text{is odd}
\end{align}
and
\begin{align}\label{1.3}
 &\left|\int_a^b f(x)dx-\frac{f(a)+f(b)}{2}(b-a)+\sum_{i=1}^{ \lfloor \frac{n-1}{2} \rfloor}\frac{2i(b-a)^{2i+1}}{2^{2i}(2i+1)!}f^{(2i)}\left(\frac{a+b}{2}\right)\right|\nonumber\\
 \leq & \frac{(b-a)^{n+1}n}{(n+1)!\,2^n}\|f^{(n)}\|_\infty,\quad
\text{if}\ n\ \text{is even}
\end{align}
where $\left\lfloor \frac{n-1}{2}\right\rfloor$ denotes  the integer part of $\frac{n-1}{2}$.

If
there exists a real number  $\gamma_n$ such that $\gamma_n\le f^{(n)}(x)$, $x\in [a, b]$, then
\begin{align}\label{1.4}
 &\left|\int_a^b f(x)dx-\frac{f(a)+f(b)}{2}(b-a)+\sum_{i=1}^{ \lfloor \frac{n-1}{2} \rfloor}\frac{2i(b-a)^{2i+1}}{2^{2i}(2i+1)!}f^{(2i)}\left(\frac{a+b}{2}\right)\right|\nonumber\\
 \leq & \left[\frac{f^{(n-1)}(b)-f^{(n-1)}(a)}{b-a}-\gamma_n\right]\frac{n-1}{n!\,2^n}(b-a)^{n+1},\quad
\text{if}\ n\ \text{is odd}.
\end{align}

If
there exists a real number  $\Gamma_n$ such that $f^{(n)}(x)\le \Gamma_n$, $x\in [a, b]$, then
\begin{align}\label{1.5}
 &\left|\int_a^b f(x)dx-\frac{f(a)+f(b)}{2}(b-a)+\sum_{i=1}^{ \lfloor \frac{n-1}{2} \rfloor}\frac{2i(b-a)^{2i+1}}{2^{2i}(2i+1)!}f^{(2i)}\left(\frac{a+b}{2}\right)\right|\nonumber\\
 \leq & \left[\Gamma_n-\frac{f^{(n-1)}(b)-f^{(n-1)}(a)}{b-a}\right]\frac{n-1}{n!\,2^n}(b-a)^{n+1},\quad
\text{if}\ n\ \text{is odd}.
\end{align}
\end{theorem}

In \cite{l2005tjm}, by  choosing
\begin{equation}\label{1.6}
 M_n(x)=\left\{ \begin{array}{ll}\displaystyle
 \frac{(x-a)^{n}}{n!},\ \ &\displaystyle x\in \left[a, \frac{a+b}{2}\right], \hfill
 \medskip\\
 \displaystyle \frac{(x-b)^{n}}{n!},\ \ &\displaystyle x\in \left(\frac{a+b}{2}, b\right], \hfill
\end{array} \right.
\end{equation}
Liu provided the following sharp perturbed midpoint inequalities
by using the variant of the Gr\"uss inequality.
\begin{theorem} \label{th1.2}
Let $f: [a, b] \rightarrow \mathbb{R}$ be such that $f^{(n)}$ is integrable with  $\gamma_n\le f^{(n)}(x) \le \Gamma_n$ for all $x\in [a, b]$, where $\gamma_n, \Gamma_n\in R$ are constants. Then
\begin{align}\label{1.7}
 &\left|\int_a^b f(x)dx-f\left(\frac{a+b}{2}\right)(b-a)-\sum_{k=1}^{n-1}\frac{[1+(-1)^k](b-a)^{k+1}}{2^{k+1}(k+1)!}f^{(k)}\left(\frac{a+b}{2}\right)\right.\nonumber\\
 &\left.
 -\frac{[1+(-1)^n](b-a)^{n+1}}{2^{n+1}(n+1)!} \frac{f^{(n-1)}(b)-f^{(n-1)}(a)}{b-a}\right|\nonumber\\
 \leq & \frac{\Gamma_n-\gamma_n}{(n+1)!}\left[\frac{1-(-1)^n}{2}+\frac{[1+(-1)^n]n}{(n+1)\sqrt[n]{n+1}}\right]\frac{1}{2^{n+1}}(b-a)^{n+1}.
\end{align}
\end{theorem}

In \cite{l2005}, by  choosing the kernel
\begin{equation}\label{1.8}
 S_n(x)=\left\{ \begin{array}{ll}\displaystyle
 \frac{(x-a)^{n-1}}{n!}\left[x-a-\frac{n(b-a)}{6}\right],\ \ &\displaystyle x\in \left[a, \frac{a+b}{2}\right], \hfill
 \medskip\\
 \displaystyle \frac{(x-b)^{n-1}}{n!}\left[x-b+\frac{n(b-a)}{6}\right],\ \ &\displaystyle x\in \left(\frac{a+b}{2}, b\right], \hfill
\end{array} \right.
\end{equation}
an inequality of Simpson type for an $n$-times continuously differentiable mapping
was given as follows.
\begin{theorem} \label{th1.3}
Let $f: [a, b] \rightarrow \mathbb{R}$  be an $n$-times continuously differentiable mapping,
$n\ge 1$ and such that $\|f^{(n)}\|_\infty:=\sup_{x\in(a,b)}|f^{(n)}(x)|<\infty$. Then
\begin{align}\label{1.9}
 &\left|\int_a^b f(x)dx-\frac{b-a}{6}\left[f(a)+4f\left(\frac{a+b}{2}\right)+f(b)\right]+\sum_{i=1}^{ \lfloor \frac{n-1}{2} \rfloor}\frac{(i-1)(b-a)^{2i+1}}{3(2i+1)!2^{2i-1}}f^{(2i)}\left(\frac{a+b}{2}\right)\right|\nonumber\\
 \leq & \frac{(b-a)^{n+1}}{(n+1)!}\|f^{(n)}\|_\infty\times
 \left\{ \begin{array}{ll}\displaystyle
   \frac{4n^n}{6^{n+1}}-\frac{n-2}{3\cdot2^n},\ \ &\displaystyle \text{if}\ n<3, \hfill
 \medskip\\
 \displaystyle \frac{n-2}{3\cdot2^n},\ \ &\displaystyle \text{if}\ n\ge 3. \hfill
\end{array} \right.
\end{align}
\end{theorem}

In \cite{pv2001}, using the well-known pre-Gr\"uss inequality, Pe$\breve{c}$ari\'c and Varo$\breve{s}$anec  obtained different error bounds for inequality \eqref{1.9}. Liu \cite{l2007} generalized inequality \eqref{1.9} and also provided an improvement of \cite{pv2001}. More recently, Shi and Liu \cite{sl2009} further derived some sharp Simpson type inequalities.

The purpose of  this paper is to give a unified generalization of some quadrature rules, which not only unify the above results about error bounds for generalized mid-point,
trapezoid and Simpson's rules, but also give some
new error bounds for other quadrature
rules as special cases. Especially, we will derive two sharp error inequalities when $n$ is an odd and an even integer, respectively.

\section{Preliminaries}

In this section we present some lemmas and notations needed in the proof of our main
results.

\begin{lemma} \label{le2.1}
Let $f: [a, b] \rightarrow \mathbb{R}$ be a function such that $f^{(n-1)}$, $n \ge 1$, is absolutely continuous. Then
\begin{align}\label{2.1}
  \int_a^b f(x)dx=&(b-a)\left[(1-\theta)f\left(\frac{a+b}{2}\right)+\theta\frac{f(a)+f(b)}{2}\right]
 \nonumber\\
  & +\sum_{i=1}^{ \lfloor \frac{n-1}{2} \rfloor}
 \frac{[1-\theta(2i+1)](b-a)^{2i+1}}{(2i+1)!2^{2i}}f^{(2i)}\left(\frac{a+b}{2}\right) +R(f),
\end{align}
for all $\theta\in [0, 1]$, where $$R(f)=(-1)^n\int_a^b G_n(x)f^{(n)}(x)dx,$$ and
\begin{equation}\label{2.2}
 G_n(x)=\left\{ \begin{array}{ll}\displaystyle
 \frac{(x-a)^{n-1}}{n!}\left[x-a-\frac{\theta n(b-a)}{2}\right],\ \ &\displaystyle x\in \left[a, \frac{a+b}{2}\right], \hfill
 \medskip\\
 \displaystyle \frac{(x-b)^{n-1}}{n!}\left[x-b+\frac{\theta n(b-a)}{2}\right],\ \ &\displaystyle x\in \left(\frac{a+b}{2}, b\right]. \hfill
\end{array} \right.
\end{equation}
\end{lemma}
\begin{proof}
We briefly sketch the proof. We introduce the notations
$$P_n(x)=\frac{(x-a)^{n-1}}{n!}\left[x-a-\frac{\theta n(b-a)}{2}\right],$$ $$ Q_n(x)=\frac{(x-b)^{n-1}}{n!}\left[x-b+\frac{\theta n(b-a)}{2}\right],$$
then one can see that $P_n$ and $Q_n$ form Appell sequences of polynomials, that is
$$P_n'(x)=P_{n-1}(x),\quad Q_n'(x)=Q_{n-1}(x),\quad P_0(x)=Q_0(x)=1. $$
Thus one can use integration by parts to prove that \eqref{2.1} holds.
\end{proof}

\begin{lemma} \label{le2.2}
The Peano kernels $G_n(t)$, satisfies

\begin{equation}\label{2.3}
 \int_a^b G_n(x) dx=\left\{ \begin{array}{ll}\displaystyle
 0,\ \ &n \ \text{odd}, \hfill
 \medskip\\
 \displaystyle  \frac{(b-a)^{n+1}}{n!\,2^n}\left(\frac{1}{n+1}-\theta\right),\ \ &\displaystyle n \ \text{even}, \hfill
\end{array} \right.
\end{equation}
\medskip
\begin{equation}\label{2.4}
 \int_a^b |G_n(x)| dx=\left\{ \begin{array}{ll}\displaystyle
 \frac{(b-a)^{n+1}}{n!\,2^n}\left(\theta-\frac{1}{n+1}\right),\ \ & \text{if}\ \theta n \ge 1, \hfill
 \medskip\\
 \displaystyle  \frac{(b-a)^{n+1}}{n(n+1)!\,2^n}\left[2(\theta n)^{n+1}-\theta n(n+1)+n\right],\ \ &\text{if}\ 1\le \theta n < 1, \hfill
\end{array} \right.
\end{equation}
\medskip
\begin{equation}\label{2.5}
 \max_{x\in [a,b]} |G_n(x)| =\left\{ \begin{array}{ll}\displaystyle
 \frac{(\theta n-1)(b-a)^{n}}{n!\,2^n},\ \ & \text{if}\ \theta n >\theta+1, \hfill
 \medskip\\
 \displaystyle  \frac{\theta^n(n-1)^{n-1}(b-a)^{n}}{n!\,2^n},\ \ &\text{if}\ 1\le \theta n \le \theta+1\ \text{and}\ n>1, \hfill
 \medskip\\
 \displaystyle  \max\{1-\theta n, \theta^n(n-1)^{n-1}\}\frac{(b-a)^{n}}{n!\,2^n},\ \ &\text{if}\  \theta n < 1\ \text{and}\ n>1, \hfill
 \medskip\\
 \displaystyle  \max\{1-\theta, \theta\}\frac{b-a}{2},\ \ &\text{if}\  n=1, \hfill
\end{array} \right.
\end{equation}
\medskip
\begin{equation}\label{2.6}
 \int_a^b G_n^2(x) dx= \frac{[\theta^2n^2(2n+1)-\theta(4n^2-1)+(2n-1)](b-a)^{2n+1}}{(2n+1)(2n-1)(n!)^22^{2n}},
\end{equation}
and
\begin{align}\label{2.7}
 &\max_{x\in [a,b]} \left|G_{2m}(x)-\frac{1}{b-a}\int_a^b G_{2m}(x) dx \right| \nonumber\\
 =&\left\{ \begin{array}{ll}\displaystyle
 \max\left\{\theta -\frac{1}{2m+1}, \left[\theta(2m-1)-\frac{2m}{2m+1}\right]\right\}\frac{(b-a)^{2m}}{(2m)!\,2^{2m}},\ \   &\text{if}\ \theta (2m-1)\ge 1, \hfill
 \medskip\\
 \displaystyle  \max\left\{\left|\theta -\frac{1}{2m+1}\right|, \left|\theta(2m-1)-\frac{2m}{2m+1}\right|,\right.
 \hfill
 \medskip\\
 \displaystyle  \quad  \quad  \quad  \quad \left.\left|\theta-\frac{1}{2m+1}-\theta^{2m}(2m-1)^{2m-1}\right|
 \right\}\frac{(b-a)^{2m}}{(2m)!\,2^{2m}},\ \   &\text{if}\ \theta (2m-1)< 1. \hfill
\end{array} \right.
\end{align}

\end{lemma}
\begin{proof}
A simple calculation gives
$$\int_a^b G_n(x) dx=  \frac{(b-a)^{n+1}}{n!\,2^{n+1}}\left(\frac{1}{n+1}-\theta\right) [1-(-1)^{n+1}],$$
from which we see that \eqref{2.3} holds.

We have
\begin{align*}
 \int_a^b |G_n(x)| dx=\int_a^{\frac{a+b}{2}} |P_n(x)| dx +\int_{\frac{a+b}{2}}^b |Q_n(x)| dx =2\int_a^{\frac{a+b}{2}} |P_n(x)| dx
 = \frac{(b-a)^{n+1}}{n!\,2^n}\int_0^1 |t^{n-1}(t-\theta n)|dt,
\end{align*}
by substitution $x=a+\frac{b-a}{2}t$. If $\theta n\ge 1$, one has
$$\int_a^b |G_n(x)| dx=\frac{(b-a)^{n+1}}{n!\,2^n}\left(\theta n\int_0^1  t^{n-1}dt-\int_0^1  t^{n}dt\right)=
\frac{(b-a)^{n+1}}{n!\,2^n}\left(\theta-\frac{1}{n+1}\right).$$
If $0\le \theta n< 1$, one gets
\begin{align*}\int_a^b |G_n(x)| dx=&\frac{(b-a)^{n+1}}{n!\,2^n}\left[\int_0^{\theta n}  t^{n-1}(\theta n-t) dt+\int_{\theta n}^1  t^{n-1}(t-\theta n) dt\right]\\=&
\frac{(b-a)^{n+1}}{n(n+1)!\,2^n}\left[2(\theta n)^{n+1}-\theta n(n+1)+n\right].\end{align*}
 By combining the above two cases, \eqref{2.4} is established.

We have
$$\max_{x\in [a,b]} |G_n(x)|=\max\left\{\max_{x\in \left[a, \frac{a+b}{2}\right]} |P_n(x)|,
 \max_{x\in \left[\frac{a+b}{2}, b\right]} |Q_n(x)|\right\}=\max_{x\in \left[a, \frac{a+b}{2}\right]} |P_n(x)|.$$
When $n=1$,
$$P_1(x)=x-a-\frac{\theta(b-a)}{2},\quad x\in \left[a, \frac{a+b}{2}\right],$$
then
$$\max_{x\in [a,b]} |G_1(x)|=\max\left\{1-\theta, \theta\right\} \frac{b-a}{2}.$$
When $n>1$, since $P_n'(x)=0$ gives $x=a$ or $x=a+\frac{\theta(n-1)(b-a)}{2}$, we divide the proof of \eqref{2.5} into three steps according to the different intervals of $\theta n$.

\textit{Case $\theta n>\theta+1$}: We have
$$a<\frac{a+b}{2}<a+\frac{\theta(n-1)(b-a)}{2}<a+\frac{\theta n(b-a)}{2}.$$
So, we can get
$$P_n(x)<0,\quad P_n(x)\ \text{is decreasing,\ \  for}\ x\in \left[a, \frac{a+b}{2}\right].$$
Therefore,
$$\max_{x\in [a,b]} |G_n(x)|=-\min_{x\in \left[a, \frac{a+b}{2}\right]} P_n(x)=-P_n\left(\frac{a+b}{2}\right)=\frac{(\theta n-1)(b-a)^{n}}{n!\,2^n}.$$

\textit{Case $1\le\theta n\le\theta+1$}: We have
$$a<a+\frac{\theta(n-1)(b-a)}{2}\le\frac{a+b}{2}\le a+\frac{\theta n(b-a)}{2}.$$
So, we can get
$$P_n(x)<0, \ \text{ for}\ x\in \left[a, \frac{a+b}{2}\right].$$
Therefore,
$$\max_{x\in [a,b]} |G_n(x)|=-\min_{x\in \left[a, \frac{a+b}{2}\right]} P_n(x)=-P_n\left(a+\frac{\theta(n-1)(b-a)}{2}\right)= \frac{\theta^n(n-1)^{n-1}(b-a)^{n}}{n!\,2^n}.$$

\textit{Case $\theta n<1$}: We have
$$a<a+\frac{\theta(n-1)(b-a)}{2}< a+\frac{\theta n(b-a)}{2}<\frac{a+b}{2}.$$
So, we can get
$$P_n(x)<0,  \ \text{ for}\   x\in \left[a, a+\frac{\theta n(b-a)}{2}\right); \quad
P_n(x)\ge 0,  \ \text{ for}\   x\in \left[a+\frac{\theta n(b-a)}{2},  \frac{a+b}{2}\right].$$
Therefore,
\begin{align*} \max_{x\in [a,b]} |G_n(x)|=&\max\left\{\left|P_n\left(\frac{a+b}{2}\right)\right|, \left|P_n\left(a+\frac{\theta(n-1)(b-a)}{2}\right)\right|\right\}
= \max\{1-\theta n, \theta^n(n-1)^{n-1}\}\frac{(b-a)^{n}}{n!\,2^n}.
\end{align*}
 By combining the above three cases, \eqref{2.5} is established.

\eqref{2.6} can be obtained by a direct calculation.

From \eqref{2.3}, we have
\begin{align*}
&\max_{x\in [a,b]} \left|G_{2m}(x)-\frac{1}{b-a}\int_a^b G_{2m}(x) dx \right|\\
=&\max_{x\in [a,b]} \left|G_{2m}(x)-\frac{(b-a)^{2m}}{(2m)!\,2^{2m}}\left(\frac{1}{2m+1}-\theta\right) \right|\\
=&\max_{x\in \left[a, \frac{a+b}{2}\right]} \left|\frac{(x-a)^{2m-1}}{(2m)!}\left[x-a-\frac{2\theta m(b-a)}{2}\right]-\frac{(b-a)^{2m}}{(2m)!\,2^{2m}}\left(\frac{1}{2m+1}-\theta\right) \right|
:= \max_{x\in \left[a, \frac{a+b}{2}\right]}|F_n(x)|
\end{align*}
and $$F_n(a)=\frac{(b-a)^{2m}}{(2m)!\,2^{2m}}\left(\theta-\frac{1}{2m+1}\right).$$
We divide the proof of \eqref{2.7} into two steps according to the different intervals of $\theta (2m-1)$.

\textit{Case $\theta (2m-1)\ge 1$}:  We have
$$F_n(a)>0,\quad a+\frac{\theta(2m-1)(b-a)}{2}\geq \frac{a+b}{2}.$$
Thus, we get
\begin{align*}&\max_{x\in \left[a, \frac{a+b}{2}\right]}|F_n(x)|=\max \left\{F_n(a), \left|F_n\left(\frac{a+b}{2}\right)\right|\right\}
=
\max\left\{\theta -\frac{1}{2m+1}, \left[\theta(2m-1)-\frac{2m}{2m+1}\right]\right\}\frac{(b-a)^{2m}}{(2m)!\,2^{2m}}.\end{align*}

\textit{Case $\theta (2m-1)< 1$}:  We have
$$a<a+\frac{\theta(2m-1)(b-a)}{2}< \frac{a+b}{2}.$$
Thus, we obtain
\begin{align*}&\max_{x\in \left[a, \frac{a+b}{2}\right]}|F_n(x)|=\max \left\{|F_n(a)|, \left|F_n\left(\frac{a+b}{2}\right)\right|,
\left|F_n\left(a+\frac{\theta(2m-1)(b-a)}{2}\right)\right|\right\}\\
=&
\max\left\{\left|\theta -\frac{1}{2m+1}\right|, \left|\theta(2m-1)-\frac{2m}{2m+1}\right|, \left|\theta-\frac{1}{2m+1}-\theta^{2m}(2m-1)^{2m-1}\right|
 \right\}\frac{(b-a)^{2m}}{(2m)!\,2^{2m}}.\end{align*}
 By combining the above two cases, \eqref{2.7} is established.
 \end{proof}

Before we end this section, we  introduce the notations
$$ I=  \int_a^b f(x)dx, $$
\begin{align*}F_n= (b-a)\left[(1-\theta)f\left(\frac{a+b}{2}\right)+\theta\frac{f(a)+f(b)}{2}\right] +\sum_{i=2}^{ \lfloor \frac{n-1}{2} \rfloor}
 \frac{[1-\theta(2i+1)](b-a)^{2i+1}}{(2i+1)!2^{2i}}f^{(2i)}\left(\frac{a+b}{2}\right).\end{align*}

\section{Main results}
We first establish three error inequalities for $f^{(n)}\in L^1[a,b]$, $L^2[a,b]$ and $L^\infty[a,b]$, respectively.

\begin{theorem} \label{th3.2}
Let $f: [a, b] \rightarrow \mathbb{R}$ be a function such that $f^{(n-1)}$, $n > 1$, is absolutely continuous on $[a,b]$. If $f^{(n)}\in L^1[a,b]$, then we have
\begin{align}\label{3.3}
 |I-F_n|
 \leq  \frac{(b-a)^{n}}{n!\,2^n}\|f^{(n)}\|_1\times
 \left\{ \begin{array}{ll}\displaystyle
 \theta n-1,\ \ & \text{if}\ \theta n >\theta+1, \hfill
 \medskip\\
 \displaystyle  \theta^n(n-1)^{n-1},\ \ &\text{if}\ 1< \theta n \le \theta+1, \hfill
 \medskip\\
 \displaystyle  \max\{1-\theta n, \theta^n(n-1)^{n-1}\},\ \ &\text{if}\  \theta n \le  1, \hfill
\end{array} \right.
\end{align}
for all $\theta\in [0, 1]$, where $\|f^{(n)}\|_1:= \int_a^b |f^{(n)}(x)| dx$ is the usual Lebesgue norm on $L^1[a, b]$.
\end{theorem}

\begin{proof}
By using the identity \eqref{2.1}, we have
\begin{align}\label{3.4}
 |I-F_n|
 = \left|\int_a^b G_n(x)f^{(n)}(x)dx\right|\leq \max_{x\in [a,b]} |G_n(x)|\int_a^b |f^{(n)}(x)|dx.
\end{align}
Consequently,  inequality \eqref{3.3} follows from \eqref{3.4} and \eqref{2.5}.
\end{proof}

\begin{remark}
If we take $\theta=1$ in  \eqref{3.3}, we get the trapezoid type inequality
\begin{align*}
 &\left|\int_a^b f(x)dx-\frac{f(a)+f(b)}{2}(b-a)+\sum_{i=1}^{ \lfloor \frac{n-1}{2} \rfloor}\frac{2i(b-a)^{2i+1}}{2^{2i}(2i+1)!}f^{(2i)}\left(\frac{a+b}{2}\right)\right|\nonumber\\
 \leq & \frac{(b-a)^{n}}{n!\,2^n}\|f^{(n)}\|_1\times
 \left\{ \begin{array}{ll}\displaystyle
    n-1,\ \ &\displaystyle \text{if}\  n> 1, \hfill
 \medskip\\
 \displaystyle 1,\ \ &\displaystyle \text{if}\ n=1. \hfill
\end{array} \right.
\end{align*}

If we take $\theta=\frac{1}{3}$ in  \eqref{3.3}, we recapture the Simpson type inequality
\begin{align*}
 &\left|\int_a^b f(x)dx-\frac{b-a}{6}\left[f(a)+4f\left(\frac{a+b}{2}\right)+f(b)\right]+\sum_{i=1}^{ \lfloor \frac{n-1}{2} \rfloor}\frac{(i-1)(b-a)^{2i+1}}{3(2i+1)!2^{2i-1}}f^{(2i)}\left(\frac{a+b}{2}\right)\right|\nonumber\\
 \leq & \|f^{(n)}\|_1\times
 \left\{ \begin{array}{ll}\displaystyle
   \frac{(n-3)(b-a)^n}{3(n!)2^n},\ \ &\displaystyle \text{if}\ n\ge 4, \hfill
 \medskip\\
 \displaystyle \frac{(b-a)^3}{324},\ \ &\displaystyle \text{if}\ n= 3, \hfill
 \medskip\\
 \displaystyle \frac{(b-a)^2}{24},\ \ &\displaystyle \text{if}\ n= 2, \hfill
 \medskip\\
 \displaystyle \frac{b-a}{3},\ \ &\displaystyle \text{if}\ n= 1, \hfill
\end{array} \right.
\end{align*}
which has been appeared in \cite[Theorem 4]{l2007}.

If we take $\theta=0$ in  \eqref{3.3}, we recapture the midpoint type inequality
\begin{align*}
   \left|\int_a^b f(x)dx-f\left(\frac{a+b}{2}\right)(b-a)-\sum_{i=1}^{ \lfloor \frac{n-1}{2} \rfloor}\frac{(b-a)^{2i+1}}{(2i+1)!2^{2i}}f^{(2i)}\left(\frac{a+b}{2}\right)\right|
 \leq   \frac{(b-a)^{n}}{n!2^{n}} \|f^{(n)}\|_1,
\end{align*}
which has been appeared in \cite[Corollary 4.15]{cd1}.

If we take $\theta=\frac{1}{2}$ in  \eqref{3.3}, we get the averaged  midpoint-trapezoid type inequality
\begin{align*}
 &\left|\int_a^b f(x)dx-\frac{b-a}{4}\left[f(a)+2f\left(\frac{a+b}{2}\right)+f(b)\right]+\sum_{i=1}^{ \lfloor \frac{n-1}{2} \rfloor}\frac{(2i-1)(b-a)^{2i+1}}{(2i+1)!2^{2i+1}}f^{(2i)}\left(\frac{a+b}{2}\right)\right|\nonumber\\
 \leq & \|f^{(n)}\|_1\times
 \left\{ \begin{array}{ll}\displaystyle
   \frac{(n-2)(b-a)^n}{n!2^{n+1}},\ \ &\displaystyle \text{if}\ n\ge 3, \hfill
 \medskip\\
 \displaystyle \frac{(b-a)^2}{32},\ \ &\displaystyle \text{if}\ n= 2, \hfill
 \medskip\\
 \displaystyle \frac{b-a}{4},\ \ &\displaystyle \text{if}\ n= 1. \hfill
\end{array} \right.
\end{align*}
\end{remark}

\begin{theorem} \label{th3.3}
Let $f: [a, b] \rightarrow \mathbb{R}$ be a function such that $f^{(n-1)}$, $n > 1$, is absolutely continuous on $[a,b]$. If $f^{(n)}\in L^2[a,b]$, then we have
\begin{align}\label{3.5}
 |I-F_n|
 \leq  \frac{(b-a)^{n+\frac{1}{2}}}{n!\,2^n}\|f^{(n)}\|_2 \sqrt{\frac{\theta^2n^2(2n+1)-\theta(4n^2-1)+(2n-1) }{(2n+1)(2n-1)}},
\end{align}
for all $\theta\in [0, 1]$, where $\|f^{(n)}\|_2:= \left(\int_a^b |f^{(n)}(x)|^2 dx\right)^{\frac{1}{2}}$ is the usual Lebesgue norm on $L^2[a, b]$.
\end{theorem}

\begin{proof}
By using the identity \eqref{2.1}, we have
\begin{align}\label{3.6}
 |I-F_n|
 = \left|\int_a^b G_n(x)f^{(n)}(x)dx\right|\leq \|f^{(n)}\|_2\|G_n\|_2.
\end{align}
Consequently,  inequality \eqref{3.5} follows from \eqref{3.6} and \eqref{2.6}.
\end{proof}

\begin{remark}
If we take $\theta=1$ in  \eqref{3.5}, we get the trapezoid type inequality
\begin{align*}
 &\left|\int_a^b f(x)dx-\frac{f(a)+f(b)}{2}(b-a)+\sum_{i=1}^{ \lfloor \frac{n-1}{2} \rfloor}\frac{2i(b-a)^{2i+1}}{2^{2i}(2i+1)!}f^{(2i)}\left(\frac{a+b}{2}\right)\right|\nonumber\\
 \leq & \frac{(b-a)^{n+\frac{1}{2}}}{n!\,2^n}\|f^{(n)}\|_2 \sqrt{\frac{n(2n^2-3n+2) }{(2n+1)(2n-1)}}.
\end{align*}

If we take $\theta=\frac{1}{3}$ in  \eqref{3.5}, we recapture the Simpson type inequality
\begin{align*}
 &\left|\int_a^b f(x)dx-\frac{b-a}{6}\left[f(a)+4f\left(\frac{a+b}{2}\right)+f(b)\right]+\sum_{i=1}^{ \lfloor \frac{n-1}{2} \rfloor}\frac{(i-1)(b-a)^{2i+1}}{3(2i+1)!2^{2i-1}}f^{(2i)}\left(\frac{a+b}{2}\right)\right|\nonumber\\
 \leq & \frac{(b-a)^{n+\frac{1}{2}}}{n!\,2^n}\|f^{(n)}\|_2 \sqrt{\frac{2n^3-11n^2+18n-6}{9(2n+1)(2n-1)}},
\end{align*}
which has been appeared in \cite[Theorem 5]{l2007}.

If we take $\theta=0$ in  \eqref{3.5}, we recapture the midpoint type inequality
\begin{align*}
   \left|\int_a^b f(x)dx-f\left(\frac{a+b}{2}\right)(b-a)-\sum_{i=1}^{ \lfloor \frac{n-1}{2} \rfloor}\frac{(b-a)^{2i+1}}{(2i+1)!2^{2i}}f^{(2i)}\left(\frac{a+b}{2}\right)\right|
 \leq   \frac{(b-a)^{n+\frac{1}{2}}}{n!\,2^n}\|f^{(n)}\|_2 \frac{1 }{\sqrt{2n+1}},
\end{align*}
which has been appeared in \cite[Corollary 4.15]{cd1}.

If we take $\theta=\frac{1}{2}$ in  \eqref{3.5}, we get the averaged  midpoint-trapezoid type inequality
\begin{align*}
 &\left|\int_a^b f(x)dx-\frac{b-a}{4}\left[f(a)+2f\left(\frac{a+b}{2}\right)+f(b)\right]+\sum_{i=1}^{ \lfloor \frac{n-1}{2} \rfloor}\frac{(2i-1)(b-a)^{2i+1}}{(2i+1)!2^{2i+1}}f^{(2i)}\left(\frac{a+b}{2}\right)\right|\nonumber\\
 \leq & \frac{(b-a)^{n+\frac{1}{2}}}{n!\,2^n}\|f^{(n)}\|_2 \sqrt{\frac{2n^3-7n^2+8n-2}{4(2n+1)(2n-1)}}.
\end{align*}
\end{remark}

\begin{theorem} \label{th3.1}
Let $f: [a, b] \rightarrow \mathbb{R}$  be an $n$-times continuously differentiable mapping,
$n\ge 1$ and such that $\|f^{(n)}\|_\infty<\infty$. Then we have
\begin{align}\label{3.1}
 |I-F_n|
 \leq  \frac{(b-a)^{n+1}}{(n+1)!2^n}\|f^{(n)}\|_\infty\times
 \left\{ \begin{array}{ll}\displaystyle
    \theta(n+1)-1,\ \ &\displaystyle \text{if}\ \theta n\ge 1, \hfill
 \medskip\\
 \displaystyle 2\theta^{n+1}n^n-\theta(n+1)+1,\ \ &\displaystyle \text{if}\ 0\le \theta n<1, \hfill
\end{array} \right.
\end{align}
for all $\theta\in [0, 1]$.
\end{theorem}

\begin{proof}
Using the identity \eqref{2.1}, we get
\begin{align}\label{3.2}
 |I-F_n|
 = \left|\int_a^b G_n(x)f^{(n)}(x)dx\right|\leq \|f^{(n)}\|_\infty\int_a^b |G_n(x)|dx.
\end{align}
Consequently, inequality \eqref{3.1} follows from \eqref{3.2} and \eqref{2.4}.
\end{proof}

\begin{remark}
If we take $\theta=1$ and $\theta=\frac{1}{3}$ in  \eqref{3.1}, respectively, we recapture the trapezoid type inequality \eqref{1.3} and the Simpson type inequality \eqref{1.9}, respectively. If we take $\theta=0$ in  \eqref{3.1}, we recapture the midpoint type inequality
appeared in \cite[Corollary 4.15]{cd1}.
\end{remark}

Next, if $f^{(n)}$ is integrable and bounded, we  prove some new error inequalities and perturbed   error inequalities, respectively.

\begin{theorem} \label{th3.3}
Let $f: [a, b] \rightarrow \mathbb{R}$ be such that $f^{(n)} (n\ge 1)$ is integrable with  $\gamma_n\le f^{(n)}(x) \le \Gamma_n$ for all $x\in [a, b]$, where $\gamma_n, \Gamma_n\in R$ are constants.

(1)\ If $n$ is an odd integer, we have
\begin{align}\label{3.7}
 |I-F_n|
 \leq  \frac{\Gamma_n-\gamma_n}{2} \frac{(b-a)^{n+1}}{(n+1)!2^n} \times
 \left\{ \begin{array}{ll}\displaystyle
    \theta(n+1)-1,\ \ &\displaystyle \text{if}\ \theta n\ge 1, \hfill
 \medskip\\
 \displaystyle 2\theta^{n+1}n^n-\theta(n+1)+1,\ \ &\displaystyle \text{if}\ 0\le \theta n<1, \hfill
\end{array} \right.
\end{align}
\begin{align}\label{3.8}
 |I-F_n|\leq  &\left[\frac{f^{(n-1)}(b)-f^{(n-1)}(a)}{b-a}-\gamma_n\right] \frac{(b-a)^{n+1}}{n!2^n}\nonumber\\
  &
 \times
 \left\{ \begin{array}{ll}\displaystyle
 \theta n-1,\ \ & \text{if}\ \theta n >\theta+1, \hfill
 \medskip\\
 \displaystyle  \theta^n(n-1)^{n-1},\ \ &\text{if}\ 1< \theta n \le \theta+1, \hfill
 \medskip\\
 \displaystyle  \max\{1-\theta n, \theta^n(n-1)^{n-1}\},\ \ &\text{if}\  \theta n \le  1, \hfill
\end{array} \right.
\end{align}
and
\begin{align}\label{3.9}
 |I-F_n|\leq  &\left[\Gamma_n-\frac{f^{(n-1)}(b)-f^{(n-1)}(a)}{b-a}\right] \frac{(b-a)^{n+1}}{n!2^n}\nonumber\\
  &
 \times
 \left\{ \begin{array}{ll}\displaystyle
 \theta n-1,\ \ & \text{if}\ \theta n >\theta+1, \hfill
 \medskip\\
 \displaystyle  \theta^n(n-1)^{n-1},\ \ &\text{if}\ 1< \theta n \le \theta+1, \hfill
 \medskip\\
 \displaystyle  \max\{1-\theta n, \theta^n(n-1)^{n-1}\},\ \ &\text{if}\  \theta n \le  1, \hfill
\end{array} \right.
\end{align}
for all $\theta\in [0, 1]$.

(2)\ If $n$ is an even integer $(n=2m)$, we have
\begin{align}\label{3.10}
 &\left|I-F_{2m}-\frac{(b-a)^{2m+1}}{(2m)!\,2^{2m}}\left(\frac{1}{2m+1}-\theta\right)\frac{f^{(2m-1)}(b)-f^{(2m-1)}(a)}{b-a}\right|\nonumber\\\leq  &\left[\frac{f^{(2m-1)}(b)-f^{(2m-1)}(a)}{b-a}-\gamma_{2m}\right] \frac{(b-a)^{2m+1}}{(2m)!\,2^{2m}}\nonumber\\
  &
 \times
 \left\{ \begin{array}{ll}\displaystyle
 \max\left\{\theta -\frac{1}{2m+1}, \left[\theta(2m-1)-\frac{2m}{2m+1}\right]\right\},\ \   &\text{if}\ \theta (2m-1)\ge 1, \hfill
 \medskip\\
 \displaystyle  \max\left\{\left|\theta -\frac{1}{2m+1}\right|, \left|\theta(2m-1)-\frac{2m}{2m+1}\right|,\right.
 \hfill
 \medskip\\
 \displaystyle  \quad  \quad  \quad  \quad \left.\left|\theta-\frac{1}{2m+1}-\theta^{2m}(2m-1)^{2m-1}\right|
 \right\},\ \   &\text{if}\ \theta (2m-1)< 1, \hfill
\end{array} \right.
\end{align}
and
\begin{align}\label{3.11}
 &\left|I-F_{2m}-\frac{(b-a)^{2m+1}}{(2m)!\,2^{2m}}\left(\frac{1}{2m+1}-\theta\right)\frac{f^{(2m-1)}(b)-f^{(2m-1)}(a)}{b-a}\right|\nonumber\\\leq  &\left[\Gamma_{2m}-\frac{f^{(2m-1)}(b)-f^{(2m-1)}(a)}{b-a}\right] \frac{(b-a)^{2m+1}}{(2m)!\,2^{2m}}\nonumber\\
  &
 \times
 \left\{ \begin{array}{ll}\displaystyle
 \max\left\{\theta -\frac{1}{2m+1}, \left[\theta(2m-1)-\frac{2m}{2m+1}\right]\right\},\ \   &\text{if}\ \theta (2m-1)\ge 1, \hfill
 \medskip\\
 \displaystyle  \max\left\{\left|\theta -\frac{1}{2m+1}\right|, \left|\theta(2m-1)-\frac{2m}{2m+1}\right|,\right.
 \hfill
 \medskip\\
 \displaystyle  \quad  \quad  \quad  \quad \left.\left|\theta-\frac{1}{2m+1}-\theta^{2m}(2m-1)^{2m-1}\right|
 \right\},\ \   &\text{if}\ \theta (2m-1)< 1, \hfill
\end{array} \right.
\end{align}
for all $\theta\in [0, 1]$.
\end{theorem}

\begin{proof}
(1)\ For $n$ odd, by \eqref{2.3} and \eqref{2.1}, we get
\begin{align}\label{3.12}
 I-F_n
 = -\int_a^b G_n(x)[f^{(n)}(x)-C]dx,
\end{align}
where $C\in R$ is a constant.

If we choose $C=\frac{\gamma_n+\Gamma_n}{2}$, we have
\begin{align}\label{3.13}
 |I-F_n|\leq \max_{x\in [a,b]}\left|f^{(n)}(x)-\frac{\gamma_n+\Gamma_n}{2}\right| \int_a^b |G_n(x)| dx= \frac{\Gamma_n-\gamma_n}{2}\int_a^b |G_n(x)| dx,
\end{align}
and hence inequality \eqref{3.7} follows from \eqref{3.13} and \eqref{2.4}.

If we choose $C=\gamma_n$, we have
\begin{align}\label{3.14}
 |I-F_n|\leq  \max_{x\in [a,b]}|G_n(x)|\int_a^b |f^{(n)}(x)-\gamma_n| dx,
\end{align}
and hence inequality \eqref{3.8} follows from \eqref{3.14} and \eqref{2.5}.

Similarly we can prove that inequality \eqref{3.9} holds.

(2)\ For $n$ even, by \eqref{2.3} and \eqref{2.1}, we can obtain
\begin{align}\label{3.15}
 &\left|I-F_{2m}-\frac{(b-a)^{2m+1}}{(2m)!\,2^{2m}}\left(\frac{1}{2m+1}-\theta\right)\frac{f^{(2m-1)}(b)-f^{(2m-1)}(a)}{b-a}\right|\nonumber\\
 =& \left|\int_a^b \left[G_{2m}(x)-\frac{1}{b-a}\int_a^b G_{2m}(x) dx\right][f^{(2m)}(x)-C]dx\right|,
\end{align}
where $C\in R$ is a constant.

If we choose $C=\gamma_{2m}$, we have
\begin{align}\label{3.16}
 &\left|I-F_{2m}-\frac{(b-a)^{2m+1}}{(2m)!\,2^{2m}}\left(\frac{1}{2m+1}-\theta\right)\frac{f^{(2m-1)}(b)-f^{(2m-1)}(a)}{b-a}\right|\nonumber\\
 \le & \max_{x\in [a,b]}\left|G_{2m}(x)-\frac{1}{b-a}\int_a^b G_{2m}(x) dx\right| \int_a^b |f^{(2m)}(x)-\gamma_{2m}|dx,
\end{align}
and hence inequality \eqref{3.10} follows from \eqref{3.16} and \eqref{2.7}.

Similarly we can prove that inequality \eqref{3.11} holds.
\end{proof}

\begin{remark}
If $n$ is an odd integer and we take $\theta=1$ in  \eqref{3.7}, \eqref{3.8} and \eqref{3.9}, respectively, we recapture inequalities \eqref{1.2}, \eqref{1.4} and \eqref{1.5}, respectively. If $n$ is an odd integer and we take $\theta=0$ in  \eqref{3.7},  we recapture inequality \eqref{1.7}.
If $n$ is an odd integer and we take $\theta=\frac{1}{3}$ in  \eqref{3.7}, \eqref{3.8} and \eqref{3.9},  respectively, we recapture  \cite[inequality (16)]{l2007},  \cite[inequality (17)]{l2007} and   \cite[inequality (18)]{l2007},   respectively. If $n$ is an even integer and we take $\theta=\frac{1}{3}$ in \eqref{3.10} and \eqref{3.11}, respectively, we recapture  \cite[inequality (19)]{l2007} and  \cite[inequality (20)]{l2007}, respectively.
For other special cases, such as $\theta=0$ or $\theta=\frac{1}{2}$, the interested reader can get some
new error bounds for other quadrature
rules, which we omit  here.
\end{remark}

Finally, we  derive two sharp error inequalities when $n$ is an odd and an even integer, respectively.

\begin{theorem} \label{th3.5}
Let $f: [a, b] \rightarrow \mathbb{R}$ be a function such that $f^{(n-1)}$  is absolutely continuous on $[a,b]$ and  $f^{(n)}\in L^2[a,b]$, where $n\ge 1$ is an odd integer. Then we have
\begin{align}\label{3.17}
 |I-F_n|
 \leq  \frac{(b-a)^{n+\frac{1}{2}}}{n!\,2^n}\sqrt{\frac{\theta^2n^2(2n+1)-\theta(4n^2-1)+(2n-1) }{(2n+1)(2n-1)}}\sqrt{\sigma(f^{(n)})},
\end{align}
for all $\theta\in [0, 1]$, where $\sigma(\cdot)$ is defined by $\sigma(f)=\|f\|_2^2-\frac{1}{b-a}\left(\int_a^b f(x) dx\right)^2. $
Inequality \eqref{3.17} is sharp in the sense that the constant
$\frac{1}{n!\,2^n}\sqrt{\frac{\theta^2n^2(2n+1)-\theta(4n^2-1)+(2n-1) }{(2n+1)(2n-1)}}$ cannot
be replaced by a smaller one.
\end{theorem}

\begin{proof}
From \eqref{2.1}, \eqref{2.3} and \eqref{2.6}, we can easily get
\begin{align*}
 &|I-F_n|
 = \left|\int_a^b G_n(x)\left[f^{(n)}(x)-\frac{1}{b-a}\int_a^b f^{(n)}(x)dx\right]dx\right|\\
 \leq & \left(\int_a^b  G_n^2(x)dx\right)^{\frac{1}{2}}
 \left(\int_a^b  \left[f^{(n)}(x)-\frac{1}{b-a}\int_a^b f^{(n)}(x)dx\right]^2dx\right)^{\frac{1}{2}}\\
 =& \left(\frac{[\theta^2n^2(2n+1)-\theta(4n^2-1)+(2n-1)](b-a)^{2n+1}}{(2n+1)(2n-1)(n!)^22^{2n}}\right)^{\frac{1}{2}}\\&\times
 \left(\|f^{(n)}\|_2^2-\frac{[f^{(n-1)}(b)-f^{(n-1)}(a)]^2}{b-a} \right)^{\frac{1}{2}}\\
 =&\frac{(b-a)^{n+\frac{1}{2}}}{n!\,2^n}\sqrt{\frac{\theta^2n^2(2n+1)-\theta(4n^2-1)+(2n-1) }{(2n+1)(2n-1)}}\sqrt{\sigma(f^{(n)})}.
\end{align*}

To prove the sharpness of \eqref{3.17}, we suppose that \eqref{3.17} holds with a constant $C>0$ as
\begin{align}\label{3.18}
 |I-F_n|
 \le C  (b-a)^{n+\frac{1}{2}}\sqrt{\sigma(f^{(n)})}.
\end{align}
We may find a function $f: [a, b] \rightarrow \mathbb{R}$ such that $f^{(n-1)}$  is absolutely continuous on $[a,b]$ as
\[
 f^{(n-1)}(x)=\left\{ \begin{array}{ll}\displaystyle
 \frac{(x-a)^{n}}{(n+1)!}\left[x-a-\frac{\theta (n+1)(b-a)}{2}\right],\ \ &\displaystyle x\in \left[a, \frac{a+b}{2}\right], \hfill
 \medskip\\
 \displaystyle \frac{(x-b)^{n}}{(n+1)!}\left[x-b+\frac{\theta (n+1)(b-a)}{2}\right],\ \ &\displaystyle x\in \left(\frac{a+b}{2}, b\right]. \hfill
\end{array} \right.
\]
It follows that
\begin{equation}\label{3.19}
 f^{(n)}(x)=G_n(x).
\end{equation}
It's easy to find that the left-hand side of inequality \eqref{3.18} becomes
\begin{equation}\label{3.20}
L.H.S. \eqref{3.18}=\frac{[\theta^2n^2(2n+1)-\theta(4n^2-1)+(2n-1)](b-a)^{2n+1}}{(2n+1)(2n-1)(n!)^22^{2n}},
\end{equation}
and the right-hand side of inequality \eqref{3.18} is
\begin{equation}\label{3.21}
R.H.S. \eqref{3.18}=\frac{1}{n!\,2^n}\sqrt{\frac{\theta^2n^2(2n+1)-\theta(4n^2-1)+(2n-1) }{(2n+1)(2n-1)}}C(b-a)^{2n+1}.
\end{equation}
It follows from \eqref{3.18}, \eqref{3.20} and \eqref{3.21} that
$$C\ge \frac{1}{n!\,2^n}\sqrt{\frac{\theta^2n^2(2n+1)-\theta(4n^2-1)+(2n-1) }{(2n+1)(2n-1)}},$$
which prove that the constant $\frac{1}{n!\,2^n}\sqrt{\frac{\theta^2n^2(2n+1)-\theta(4n^2-1)+(2n-1) }{(2n+1)(2n-1)}}$ is the best possible in \eqref{3.17}.
\end{proof}

\begin{theorem} \label{th3.6}
Let $f: [a, b] \rightarrow \mathbb{R}$ be a function such that $f^{(n-1)}$  is absolutely continuous on $[a,b]$ and  $f^{(n)}\in L^2[a,b]$, where $n>1$ is an even integer $(n=2m)$. Then we have
\begin{align}\label{3.22}
 &\left|I-F_{2m}-\frac{(b-a)^{2m+1}}{(2m)!\,2^{2m}}\left(\frac{1}{2m+1}-\theta\right)\frac{f^{(2m-1)}(b)-f^{(2m-1)}(a)}{b-a}\right|\nonumber\\
 \leq  &\frac{(b-a)^{2m+\frac{1}{2}}}{(2m)!\,2^{2m}}\sqrt{\frac{[\theta^2(2m)^2(4m+1)-\theta(16m^2-1)+(4m-1)]-
 (16m^2-1)(\frac{1}{2m+1}-\theta)^2 }{(4m+1)(4m-1)}}\sqrt{\sigma(f^{(2m)})},
\end{align}
for all $\theta\in [0, 1]$. Inequality \eqref{3.22} is sharp in the sense that the constant
$$\frac{1}{(2m)!\,2^{2m}}\sqrt{\frac{[\theta^2(2m)^2(4m+1)-\theta(16m^2-1)+(4m-1)]-
 (16m^2-1)(\frac{1}{2m+1}-\theta)^2 }{(4m+1)(4m-1)}}$$ cannot
be replaced by a smaller one.
\end{theorem}

\begin{proof}
From \eqref{2.1}, \eqref{2.3} and \eqref{2.6}, we can easily get
\begin{align*}
 &\left|I-F_{2m}-\frac{(b-a)^{2m+1}}{(2m)!\,2^{2m}}\left(\frac{1}{2m+1}-\theta\right)\frac{f^{(2m-1)}(b)-f^{(2m-1)}(a)}{b-a}\right|\\
 = &\left|\int_a^b G_{2m}(x)f^{(2m)}(x)dx-\frac{1}{b-a}\int_a^b G_{2m}(x)dx\int_a^b f^{(2m)}(x)dx\right|\\
  = &\frac{1}{2(b-a)}\left|\int_a^b \int_a^b [G_{2m}(x)-G_{2m}(t)][f^{(2m)}(x)-f^{(2m)}(t)]dxdt\right|\\
 \leq & \frac{1}{2(b-a)}\left(\int_a^b \int_a^b [G_{2m}(x)-G_{2m}(t)]^2dxdt\right)^{\frac{1}{2}}
 \left(\int_a^b \int_a^b [f^{(2m)}(x)-f^{(2m)}(t)]^2dxdt\right)^{\frac{1}{2}}\\
  = &  \left( \int_a^b G_{2m}^2(x)dx-\frac{1}{b-a}\left[\int_a^bG_{2m}(t)dt\right]^2\right)^{\frac{1}{2}}
 \left( \int_a^b [f^{(2m)}(x)]^2dx-\frac{1}{b-a}\left[\int_a^bf^{(2m)}(t)dt\right]^2\right)^{\frac{1}{2}}\\
  =&\frac{(b-a)^{2m+\frac{1}{2}}}{(2m)!\,2^{2m}}\sqrt{\frac{[\theta^2(2m)^2(4m+1)-\theta(16m^2-1)+(4m-1)]-
 (16m^2-1)(\frac{1}{2m+1}-\theta)^2 }{(4m+1)(4m-1)}}\sqrt{\sigma(f^{(2m)})}.
\end{align*}

To prove the sharpness of \eqref{3.22}, we suppose that \eqref{3.22} holds with a constant $C>0$ as
\begin{align}\label{3.23}
 &\left|I-F_{2m}-\frac{(b-a)^{2m+1}}{(2m)!\,2^{2m}}\left(\frac{1}{2m+1}-\theta\right)\frac{f^{(2m-1)}(b)-f^{(2m-1)}(a)}{b-a}\right|\nonumber\\
 \le &C  (b-a)^{2m+\frac{1}{2}}\sqrt{\sigma(f^{(2m)})}.
\end{align}
We may find a function $f: [a, b] \rightarrow R$ such that $f^{(n-1)}$  is absolutely continuous on $[a,b]$ as
\begin{align*}
 &f^{(2m-1)}(x)\\=&\left\{ \begin{array}{ll}\displaystyle
 \frac{(x-a)^{2m}}{(2m+1)!}\left[x-a-\frac{\theta (2m+1)(b-a)}{2}\right] -\frac{(b-a)^{2m+1}}{(2m)!\,2^{2m+1}}\left(\frac{1}{2m+1}-\theta\right),\ \ &\displaystyle x\in \left[a, \frac{a+b}{2}\right], \hfill
 \medskip\\
 \displaystyle \frac{(x-b)^{2m}}{(2m+1)!}\left[x-b+\frac{\theta (2m+1)(b-a)}{2}\right]+
 \frac{(b-a)^{2m+1}}{(2m)!\,2^{2m+1}}\left(\frac{1}{2m+1}-\theta\right),\ \ &\displaystyle x\in \left(\frac{a+b}{2}, b\right]. \hfill
\end{array} \right.
\end{align*}
It follows that
\begin{equation}\label{3.24}
 f^{(2m)}(x)=G_{2m}(x).
\end{equation}
It's easy to find that the left-hand side of inequality \eqref{3.23} becomes
\begin{align}\label{3.25}
&L.H.S. \eqref{3.23}\nonumber\\=&\frac{[\theta^2(2m)^2(4m+1)-\theta(16m^2-1)+(4m-1)]-
 (16m^2-1)(\frac{1}{2m+1}-\theta)^2 }{(4m+1)(4m-1)((2m)!)^22^{4m}}(b-a)^{4m+1},
\end{align}
and the right-hand side of inequality \eqref{3.23} is
\begin{align}\label{3.26}
&R.H.S. \eqref{3.23}\nonumber\\=&\frac{1}{(2m)!\,2^{2m}}\sqrt{\frac{[\theta^2(2m)^2(4m+1)-\theta(16m^2-1)+(4m-1)]-
 (16m^2-1)(\frac{1}{2m+1}-\theta)^2 }{(4m+1)(4m-1)}}C(b-a)^{4m+1}.
\end{align}
It follows from \eqref{3.23}, \eqref{3.25} and \eqref{3.26} that
$$C\ge \frac{1}{(2m)!\,2^{2m}}\sqrt{\frac{[\theta^2(2m)^2(4m+1)-\theta(16m^2-1)+(4m-1)]-
 (16m^2-1)(\frac{1}{2m+1}-\theta)^2 }{(4m+1)(4m-1)}},$$
which prove that the constant $\frac{1}{(2m)!\,2^{2m}}\sqrt{\frac{[\theta^2(2m)^2(4m+1)-\theta(16m^2-1)+(4m-1)]-
 (16m^2-1)(\frac{1}{2m+1}-\theta)^2 }{(4m+1)(4m-1)}}$ is the best possible in \eqref{3.22}.
\end{proof}

\begin{remark}
If  we take $\theta=\frac{1}{3}$ in  \eqref{3.17} and \eqref{3.22}, respectively, we recapture sharp Simpson type inequalities \cite[inequality (24)]{sl2009} and \cite[inequality (29)]{sl2009}, respectively.
For other special cases, such as $\theta=0$, $\theta=1$ or $\theta=\frac{1}{2}$, the interested reader can get some
new sharp error inequalities for other quadrature
rules, which we omit  here.
\end{remark}

\section*{Acknowledgements}

This work was supported by the National Natural Science Foundation of China (Grant No. 41174165, 40975002), the Tianyuan Fund
of Mathematics (Grant No. 11026211) and the Natural Science Foundation of the Jiangsu Higher Education Institutions (Grant No. 09KJB110005).

\end{document}